\documentclass{amsart}
\usepackage[all]{xy}
\usepackage{amsmath,amssymb,amsthm}
\usepackage{a4,graphicx}
\begin{document}
\newcommand{\comment}[1]{}
\newtheorem{theorem}{Theorem}[section]
\newtheorem{lemma}[theorem]{Lemma}
\newtheorem{remark}[theorem]{Remark}
\newtheorem{definition}[theorem]{Definition}
\newtheorem{corollary}[theorem]{Corollary}
\newtheorem{example}[theorem]{Example}
\newtheorem{proposition}[theorem]{Proposition}

\font\MBOLD=cmmib11

\title[Isospectral surfaces]
{Isospectral surfaces of genus 2 and 3}

\author{Hyunsuk Kang}
\address{Department of Pure Mathematics and Mathematical Statistics, 
University of Cambridge, Cambridge, UK}
\email{hk244@dpmms.cam.ac.uk}

\date{Dec 20, 2005}

\begin{abstract}
We give examples of isospectral non-isometric surfaces of genus 2 and 3 with variable curvatures and apply the result to construct isospectral potentials 
on Riemann surfaces of genus 2.
\end{abstract}

\maketitle

\section{Introduction}

In this paper, we consider compact isospectral Riemannian 2-manifolds without boundary and we reserve the term Riemann surfaces for ones with constant 
curvature -1.  Many authors, including  Brooks-Tse \cite{refBT}, Buser \cite{refBu1} and Vign\'eras \cite{refVi} have constructed pairs of isospectral 
non-isometric Riemann surfaces of genus $\ge$ 4 as well as surfaces of genus 3 with variable negative curvature.  
Such constructions have developed rapidly since Sunada's theorem \cite{refSu} was proved realising that certain group theoretic properties provide 
isospectrality of manifolds.  
The main theorem of this paper is the following.

\begin{theorem}\label{thmg2}
There exists a pair of isospectral non-isometric surfaces of genus 2 with variable negative curvature.
\end{theorem}

We also give new examples of such surfaces of genus 3 in the section \ref{genus3}.  We refer \cite{refBu} for preliminary theorems and their proofs.  
First, Sunada's theorem, the main tool of our construction, is discussed in the next section.

\section{Sunada's Theorem}\label{secsu}

Sunada's original proof is based on properties of zeta functions in number field theory with the fact that the same zeta functions give the same eigenvalue 
spectrum.  There is more than one version of the theorem and we approach it from the group theoretic point of view with the following definition 
concerning finite groups.

\begin{definition}\label{almost}
Let T be a finite group, and let $U$, $V$ be subgroups of T. Then
\begin{enumerate}
\item[$(i)$]
$U$ and $V$ are almost conjugate or Gassmann equivalent if
\begin{equation*}
|[g] \cap U| = |[g] \cap V|,
\end{equation*}
where $[g]$ denotes the conjugacy class of g in T and $|\cdot|$ denotes the cardinality of a set, in which case, $U$ and $V$ are called Gassmann subgroups 
and
\item[$(ii)$]
$(T, U, V)$ is called a Sunada triple if in addition, $U$ and $V$ are non-conjugate in T.
\end{enumerate}
\end{definition}

The groups in a Sunada triple are arranged to be covering transformation groups of normal coverings and they are related with heat kernels for the 
corresponding covering manifolds to induce the trace formula.  Since the heat trace determines the spectrum, manifolds with the same heat trace are 
isospectral.

\begin{theorem}[Sunada]\label{Sunada}
Let $U$ and $V$ be almost conjugate subgroups in a finite group T, and let $\phi:\pi_{1}(M_0) \longrightarrow T$ be a surjective homomorphism 
where $M_0$ is a compact manifold.  If $M_{1}$ and $M_{2}$ are the coverings of $M_0$ with $\pi_{1}(M_1) \cong \phi^{-1}(U)$ and 
$\pi_{1}(M_2) \cong \phi^{-1}(V)$,
then $\zeta_{M_{1}}(s)$ = $\zeta_{M_{2}}(s)$.  In particular, $M_{1}$ is isospectral to $M_{2}$.
\end{theorem}

The following diagram represents the stated coverings.

\begin{displaymath}\label{cover1}
\xymatrix{ &  M \ar[dl]_U  \ar[dd]^T   \ar[dr]^V & \\
M_1  \ar[dr]  & & \ar[dl]  M_2\\
& M_0 & }
\end{displaymath}

Note that if $U$ and $V$ are conjugate, then $M_{1}$ is isometric to $M_{2}$.  Hence to construct a pair of isospectral non-isometric surfaces, it is a
necessary condition that $H_1$ and $H_2$ are non-conjugate.  Even then, it is not sufficient enough to ensure a non-isometry between the covering manifolds 
$M_1$ and $M_2$.  One way of checking the existence of non-isometry between manifolds is by comparing the behaviour of geodesics on them and this is 
certainly the most common method used in Riemann surfaces since Huber's theorem states that compact Riemann surfaces have the same eigenvalue spectrum 
if and only if they have the same geodesic length spectrum.

For a non-isometry, Sunada imposed `bumpy metrics', that is, nowhere locally homegeneous metrics, on $M_0$.  Sunada \cite{refSu} showed that given any 
metric, there exists a nowhere locally homogeneous metric arbitrarily close in the $C^\infty$-topology to the given metric.

To complete Sunada's construction, we need a surjective homomorphism from $\pi_{1}(M_0)$ onto $T$, the group of isometries of the covering manifold $M$ in 
the diagram.  This is an easy task for 2 dimensional manifolds but the case for dimension 3 requires more care.  
For dimension 4 and higher, one can easily obtain 
the required homomorphism by adding handles and apply Sunada's theorem to have isospectral manifolds.   
Note that Sunada's theorem can be extended to orbifolds.

\section{Isospectral surfaces of small genus}

As mentioned in the beginning, Sunada's method was applied to construct isospectral Riemann surfaces of genus 4 and above in \cite{refBT} 
and \cite{refBu}.  In this section, we explain how we obtain the Sunada triples with desired genera 2 and 3.  The difficulty of having such surfaces is that 
there are only finitely many possibilities for orders of generators of $T$ in Definition \ref{almost}.  Moreover, among the Sunada triples satisfying 
these conditions, the requirement of the smoothness of $M_i$, $i=1,2$ restricts the choices.

Consider the coverings of $M_0$ with covering transformation groups $U$, $V$, and $T$ as in the diagram above.  
We give a brief summary of the constructions of $M$ and $M_i$, $i=0,1,2$ and for details, see Buser, Ch.12 \cite{refBu}.

\subsection{Construction of $M_0$}

Consider a geodesic polygon $\mathcal{D}$ in the hyperbolic plane with $N$ pairs of edges $(\alpha_i, \alpha'_i)$ where edges of each pair are of 
the same length.  We parameterise each edge by the unit interval with constant speed $a_{i}:[0,1] \longrightarrow \alpha_i$.  
Let $\mathcal{R}$ be the identification $a_i(t)=a'_i(1-t)$ by an orientation reversing isometry for each pair, 
and let $M_0$ be the quotient $\mathcal{D}/ \mathcal{R}$.  In general, $M_0$ is a topological surface with singularities.

The gluing process assigns an equivalence class $[p]:=(p_1,..., p_k)$ of vertices of $\mathcal{D}$ called a vertex cycle to an vertex $p^*$ of $M_0$ and we 
associate an element $g_{[p]}$ of $T$ with each vertex cycle $[p]$.  If the sum of angles at $[p]$ is $2 \pi / m$ where $m \in \mathbb{Z}^+$, 
then $p^*$ is said to have the order $m$ and denote its order by $ord[p]$.

\subsection{Construction of $M_1$}

Let $U$ be a subgroup of $T=<A_{1}, ..., A_{n}>$ with the index $[T:U]=k$.  Let $g_1$=$id$,..., $g_k$ be the coset representatives for $U$ in $T$ and we 
take one copy $\mathcal{D}_{g_i}$ for each $g_i$.  Identify the same type of edges in $\mathcal{D}_{g_i}$ and $\mathcal{D}_{g_j}$ whenever 
$U_{g_j}=U_{g_i A_m}$.  Then $M_1$ is formed as the quotient of $\sqcup_{i=1}^{k}\mathcal{D}_{g_i}$ by this identification and 
$M_2$ is similarly constructed.

Suppose that a vertex cycle $[p]$ yields to $p^*$ in $M_0$ and its lift $\tilde{p}$ in $M$ is projected to $\hat{p}$ in $M_1$.  In order to have a smooth 
surface $M_1$, all isotropy groups are required to be identities, otherwise the resulting surfaces are orbifolds.  Since $T$ and $U$ act freely on 
$\sqcup_{g \in T}\mathcal{D}_{g}$, the only possible non-trivial isotropy groups are vertex stabilisers and we have the following lemma.

\begin{lemma}\label{tr}
Let T be the covering transformation group of $p:M \rightarrow M_0$ and suppose that a hyperbolic polygon for $M_0$ has boundary edges identified 
in N pairs and mapped to generators of T.
\begin{enumerate}
\item[$(i)$] $M_1=M/U$ is smooth at the points $\hat{p}$ over $p^*$  if, and only if,
\begin{displaymath}
g_{[p]}=id, \quad \textrm{or} \quad  [(g_{[p]})^r] \cap U=\emptyset  \quad  \textrm{for $r=1,..., ord[p]-1$}.
\end{displaymath}
\item[$(ii)$] Let $U$ be a subgroup of $T$.  Then
\begin{equation*}
\chi(M_1)=[T:U](1-N+\sum_{[p]} \frac{1}{ord[p]}).
\end{equation*}
\end{enumerate}
\end{lemma}

\begin{figure}[!hbtp]
\begin{center}
\scalebox{0.65}
{
\includegraphics[width=8cm]{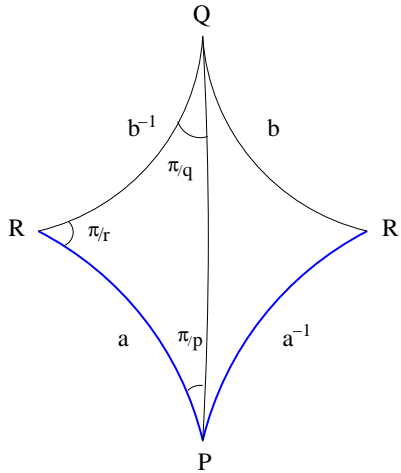}
}
\caption{A sphere with three singularities}
\label{4gon}
\end{center}
\end{figure}

Figure \ref{4gon} is an example of $M_0$ and as in \cite{refBT}, let $T$ be the index 2 subgroup of the group generated by reflections in the sides of a 
hyperbolic triangle with angles $\frac{\pi}{p}$ and $\frac{\pi}{q}$, and $\frac{\pi}{r}$.  Then $T$ is generated by $a$ and $b$ of order $p$ and $q$, 
respectively and the edges correspend to $a, b$ and their inverses so that $g_{P}=a$, $g_{Q}=b$ and $g_{R}=(ab)^{-1}$.  We use this geodesic polygon for the 
surfaces of genus 2.

\subsection{Isospectral surfaces of genus 2}\label{genus2}

In \cite{refBS}, Bosma and de Smit proved that there are exactly 19 Sunada triples of at most index 15 up to isomorphism and among their list, 
$G_{2,2,3} \cap A_{12}$ is the group which we consider here.  Note that the group $G_{2,2,3}$ is deduced from the general construction of the 
3-step abelian group $G_{p,m,k}$ in \cite{refDL}.

Let $T$ be a subgroup of $A_{12}$ with the generators
\begin{eqnarray*}
a &=& (0,7,11)(1,5,6)(2,9,10)(3,4,8),\\
b &=& (0,4,2)(1,5,9)(3,7,11)(6,10,8),\\
c &=& (ab)^{-1}=(0, 10, 5, 6, 4, 11)(1, 2, 3, 7, 8, 9).
\end{eqnarray*}

Then $T$ is the group $G_{2,2,3} \cap A_{12}$ of order 96.

In order to have smooth surfaces, the subgroups are required to avoid the conjugacy classes where the generators $a,b$ and $c$ and 
their powers belong to.  There exist such subgroups and we choose $U$ and $V$ to be

\begin{eqnarray*}
U &=&  \{ Id, (1, 7)(4, 10), (2, 8)(5, 11), (1, 7)(2, 8)(4, 10)(5, 11),\\
&& (0, 9)(2, 11, 8, 5)(3, 6)(4, 10), (0, 9)(1, 7)(2, 11, 8, 5)(3, 6),\\
&& (0, 9)(2, 5, 8, 11)(3, 6)(4, 10), (0, 9)(1, 7)(2, 5, 8, 11)(3, 6) \},\\\\
V &=& \{ Id, (1, 7)(4, 10), (0, 6)(3, 9), (0, 6)(1, 7)(3, 9)(4, 10),\\
&&  (0, 3, 6, 9)(1, 10)(2, 8)(4, 7), (0, 3, 6, 9)(1, 4)(2, 8)(7, 10),\\
&&  (0, 9, 6, 3)(1, 10)(2, 8)(4, 7), (0, 9, 6, 3)(1, 4)(2, 8)(7, 10) \}.
\end{eqnarray*}

Clearly $U$ and $V$ are of index 12.  To construct desired surfaces, we use the geodeic polygon in Figure \ref{4gon} as the fundamental domain 
whose vertices correspond to $a,b$ and $(ab)^{-1}$.  From the conjugacy classes of $T$, one can check that they are almost conjugate but not conjugate.  
Also no conjugates of $a,b$ and $(ab)^{-1}$ and their powers lie in $U$ and $V$ so that $M_1=M/U$ and $M_2=M/V$ are smooth surfaces.  
Therefore $(T, U, V)$ is a Sunada triple of index 12.  By Lemma \ref{tr} (ii),
\begin{equation*}
\chi(M_i)=12(1-2+\frac{1}{3}+\frac{1}{3}+\frac{1}{6})= -2  \quad \textrm{for $i=1,2$},
\end{equation*}
and by Sunada's theorem, they are isospctral Riemann surfaces of genus 2.

Now we investigate whether they are isometric. Studying graphs, Schreier graphs in particular, gives useful informations to determine whether surfaces 
arisen by Sunada triple are isometric.  
If there is a reflection symmetry and an orientation reversing isometry in the fundamental domain, then the surfaces are isometric whenever the 
corresponding Schreier graphs are isomorphic.  If there is no such isometry, then geodesics are taken into consideration.  For $T$, the only possibility 
for the fundamental domain is the sphere with three singularities since it has essentially two generators $a,b$.

\begin{figure}[!htbp]
\begin{center}
\includegraphics[width=8cm]{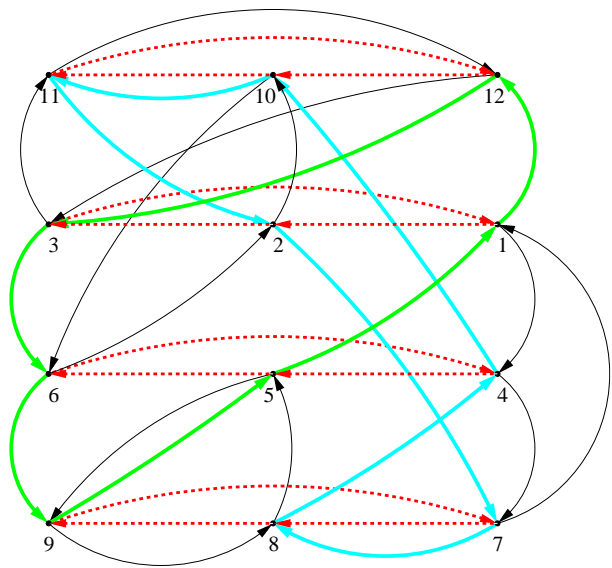}
\caption{Schreier Graph for $M_1$}
\label{graph3}
\end{center}
\end{figure}

\begin{figure}[!htbp]
\begin{center}
\includegraphics[width=8cm]{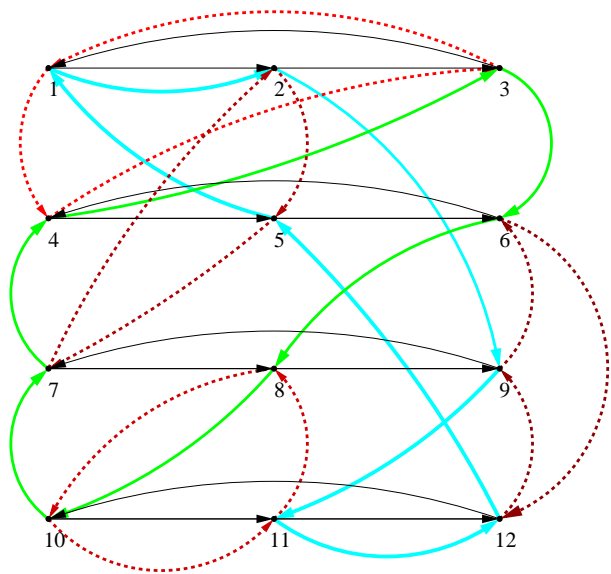}
\caption{Schreier Graph for $M_2$}
\label{graph1}
\end{center}
\end{figure}

The Schreier graphs for $M_1$ and $M_2$ are Figure \ref{graph3} and Figure \ref{graph1}, respectively and 
thin solid lines, dotted lines and thick shaded lines represent the actions of $a$, $b$ and $c$, respectively.  
The vertices in Figure \ref{graph3} are arranged so that one can see that two graphs are isomorphic with only 
arrows reversed.  Since Figure \ref{4gon} has an orientation reversing isometry, the resulting Riemann surfaces are isometric.  
To complete the proof of Theorem \ref{thmg2}, we give a nowhere locally homogeneous metric on $M_0$ so that the isospectral surfaces $M_1$ and $M_2$ 
are non-isometric with variable negative curvature.

If $M_i$, $i=1,2$, are not required to be smooth, then even the well known group $H:=\mathbb{Z}/8^{\times} \ltimes \mathbb{Z}/8$ 
gives isospectral orbifolds of genus 2.

Let  $U_{1}=\{(1,0), (3, 0), (5, 0), (7, 0)\}$ and $U_{2}=\{(1,0), (3, 4), (5, 4), (7, 0)\}$. 
Then $(H, U_{1}, U_{2})$ is a Sunada triple and we take a=(3, 2), b=(7, 1) and c=(7, 2) as the generators of $H$.  We use the sphere with four singularities 
corresponding to the vertex cycles $a, b, c$ and $abc$ of order 2,2,2 and 4 respectively.  By counting the number of intersections bewteen the conjugacy 
classes of powers of $a, b, c, abc$ and $U_{i}$, one can check that they are indeed a Sunada triple.  
Moreover $[(3,2)] \cap U_1=(3, 0)$, $[(3,2)] \cap U_2=(3, 4)$ and $[(7,2)] \cap U_1=(7, 0)$.  
Therefore $M_i$, $i=1,2$, are isospectral orbifolds genus 2 with two singularities of order two.

\subsection{Isospectral surfaces of genus 3}\label{genus3}

Brooks and Tse \cite{refBT} constructed a pair of isospectral surfaces of genus 3 with bumpy metrics on them by considering $SL(3,\mathbb{Z}/2)$.
We present an example which gives this property with a different group, namely $GL(2,\mathbb{Z}/4)$ as listed in \cite{refBS}.

Let $a$ and $b$ and their product $c$ be

\begin{displaymath}
a=
\left(\begin{array}{cc}
3 & 2  \\3 & 3
\end{array}\right), \quad
b=
\left(\begin{array}{cc}
1 & 3  \\2 & 3
\end{array}\right), \quad
c=
\left(\begin{array}{cc}
3 & 3  \\1 & 2
\end{array}\right).
\end{displaymath}

Then $a, b$ and $c$ have order 4, 4 and 6, respectively and they generate $GL(2,\mathbb{Z}/4)$.  For this group, there is a Gassmann pair $(U_{1}, U_{2})$ 
of order 8 which do not meet the conjugacy classes where the powers of $a, b$ and $c$ belong to.

\begin{eqnarray*}
U_{1} &=& \langle
\left(\begin{array}{cc}
1 & 2  \\0 & 3
\end{array}\right)
\left(\begin{array}{cc}
1 & 1  \\0 & 3
\end{array}\right)
\rangle\\
U_{2} &=& \langle
\left(\begin{array}{cc}
1 & 0  \\2 & 3
\end{array}\right)
\left(\begin{array}{cc}
1 & 0  \\1 & 3
\end{array}\right)
\rangle
= U_{1}^{T}
\end{eqnarray*}

Let the coset representatives be $<g_1, g_2>$ of order 12 where

\begin{displaymath}
g_1:=
\left(\begin{array}{cc}
3 & 3  \\1 & 2
\end{array}\right), \quad \textrm{and} \quad
g_2:=
\left(\begin{array}{cc}
3 & 0  \\0 & 3
\end{array}\right).
\end{displaymath}

Again we use the sphere with three singularities of order 4, 4 and 6 and the symmetry of the fundamental domain forces the isometry between $M/U_{1}$ and 
$M/U_{2}$.  Hence a bumpy metric on $M_0$ gives non-isometric isospectral surfaces of genus 3 with variable curvature.

\section{Isospectral potentials on a surface of genus 2}

Two smooth functions $\phi_1$ and $\phi_2$ on a smooth manifold $M$ are said to be
isospectral potentials if the spectrum of $\Delta+ \phi_1$ is identical to that of $\Delta+ \phi_2$.  
Guillemin and Kazdan \cite{refGK} showed that for compact negatively curved 2-manifolds
with simple length spectrum, the spectrum of $\Delta + \phi$
uniquely determine the potential $\phi$ and later Croke and
Sharafutdinov \cite{refCrSh} showed that it is also true for such
manifolds of any dimension.  But it was pointed out by Brooks
\cite{refBr2} that the simple length spectrum condition is never
satisfied by a metric of constant negative curvature and hence the spectral
rigidity of potential does not hold for Riemann surfaces.
He used this fact to give isospectral potentials on a surface of genus 3 and we construct ones for genus 2.

Let the Sunada triple for genus 2 be $(T, U, V)$ in the section
\ref{genus2} and let $M_1$ and $M_2$ be the normal coverings of $M_0$.  Then we showed that there exists an orientation reversing
isometry $\psi$ between $M_1$ and $M_2$.

As noticed in \cite{refBr1}, Sunada's theorem can also be applied to a Schrodinger operator $\Delta + \phi$ where $\phi$ is a smooth function on $M_0$.

\begin{theorem}[\cite{refBr1}]
For a Sunada triple $(T, U, V)$, the spectra of $\Delta + \phi_i$ on $M_i$ for $i=1,2$, coincide.  
\end{theorem}  

For the proof, one simply replaces the heat kernels of the corresponding Schrodinger operators for the heat kernels of the Laplacian.

\begin{displaymath}\label{potential}
\xymatrix{
M_1  \ar[dr]_{\pi_1} \ar[rr]^{\psi} && \ar[dl]^{\pi_2}  M_2\\
& M_0 & }
\end{displaymath}

\begin{theorem}\label{potential2}
There exist distinct isospectral potentials on a Riemann surface of genus 2.
\end{theorem}

\begin{proof}
Choose a smooth function $\phi$ on $M_0$ which does not commute with the isometry $\psi$, that is, $\phi_1 \ne \phi_2 \circ \psi$
where $\phi_i$ is the lift of $\phi$ on $M_i$, $i=1,2$.  Thus $\phi$ is not invariant under orientation reversing isometry.  Since 
the spectra of $\Delta + \phi_i$ on $M_i$ for $i=1,2$, coincide, $\phi_1$ and $\phi_2 \circ \psi$ are isospectral potentials on $M_1$ and they are not 
identical.
\end{proof}

\end{document}